\documentstyle[12pt]{article}

\newcommand{\eop}{\hfill $\#$\vspace{.1in}}

\newtheorem{property}{Property}

\newtheorem{lemma}{Lemma}

\newtheorem{conjecture}[lemma]{Conjecture}

\newcommand{\zz}{{\bf Z}}

\newcommand{\Mm}{{\cal M}}

\newcommand{\Ll}{{\cal L}}

\begin{document}

\section*{Some properties of the theory of $n$-categories}

Carlos Simpson
\newline
CNRS, Universit\'e de Nice-Sophia Antipolis

Much interest has recently focused on the problem of comparing different
definitions of $n$-categories.
\footnote{We always mean ``weak $n$-category''.}
Leinster has made a useful compendium of 10 definitions \cite{Leinster};
May has proposed another definition destined among other things to make
comparison easier \cite{May}, and he has also led the creation of an umbrella
research group including comparison as one of the main research
topics; Batanin has started a comparison of his definition
with Penon's \cite{Batanin}; and Berger has made a comparison between
Batanin's theory and the homotopy theory of spaces, introducing
techniques which should allow for other comparisons starting with
Batanin's theory \cite{Berger}. The comparison
question was first explicitly mentionned by Grothendieck in
\cite{PursuingStacks}, in a prescient prediction that many different
people would come up with different definitions of $n$-category; and
this theme was again brought up by Baez and Dolan in
\cite{Categorification}.  

The purpose of this short note is to make some observations about
properties
which one can expect any theory of $n$-categories to have, and to
conjecture that these properties characterize the theory of
$n$-categories. As a small amount of evidence for this conjecture, we
show how to go from these properties to the composition law between
mapping objects in an $n$-category. 

While we don't give the proofs here, it is not too hard to see that
Tamsamani's definition of $n$-category satisfies the properties listed
below.  
We conjecture that the other definitions satisfy these properties too.
This conjecture plus the conjecture of the previous paragraph would give
an answer to the comparison question.

Rather than giving all of the references for the various different
definitions of $n$-category, we refer the reader to
Leinster's excellent bibliography \cite{Leinster}.

The fundamental tool which we use is the Dwyer-Kan localization
\cite{DwyerKan}. 
This is a generalization of the classical Gabriel-Zisman localization,
which keeps higher homotopy data. Dwyer and Kan obtain a mapping space
between two objects, where Gabriel and Zisman obtain only the set of
homotopy classes of maps i.e. the $\pi _0$ of the mapping space. The
fundamental observation of Dwyer and Kan is that the mapping spaces
(plus their composition and higher homotopy coherence information) are
determined by the data of the original category and the localizing
subcategory.

We apply this observation to the problem of comparing definitions of
$n$-categories. In any reasonable theory of $n$-categories, one starts
with a collection of morphisms between objects which are the ``strict''
morphisms or morphisms strictly respecting the structure; these are of
course not expected to be invariant in any way from one theory to the
other. They form a usual $1$-category $nCat$.
Some of these morphisms are ``$n$-equivalences'', the $n$-categorical
generalization of equivalences of categories. 

Tamsamani in his thesis looked at $Ho-nCat$, the Gabriel-Zisman
localization dividing the $1$-category $nCat$ out by the $n$-equivalences. 
Following the philosophy of Dwyer-Kan, we look instead at $\Ll _n$ which
is the simplicial category obtained by Dwyer-Kan localization dividing
the same $nCat$ by the same subcategory of strict $n$-equivalences. Note that the
$1$-truncation of this simplicial category (corresponding to replacing
the morphism spaces by their $\pi _0$) is just the Gabriel-Zisman
localization $Ho-nCat$.  

Baez and Dolan \cite{Categorification} gave a first
precise statement of the comparison problem by saying that in any 
theory ${\bf T}$ of $n$-categories, one will have an underlying $1$-category
$nCat_{\bf T}$
and a notion of (strict) $n$-equivalence in there. They say that one should ask
whether the Gabriel-Zisman localizations for two theories, are the same.
Again, following the Dwyer-Kan philosophy, we can reinforce this question
as follows: let $\Ll _{n,{\bf T}}$ denote the Dwyer-Kan simplicial
localization
obtained by dividing $nCat_{\bf T}$ by the strict $n$-equivalences. A
{\em comparison} between theories ${\bf T}$ and ${\bf T}'$ is then an
equivalence of simplicial categories
\footnote{To be precise, an ``equivalence of simplicial categories'' can mean an
isomorphism in the Gabriel-Zisman localization obtained by dividing out
the category of simplicial categories by their equivalences.}
$$
\Ll _{n,{\bf T}} \cong \Ll _{n,{\bf T}'}.
$$

It would be possible to do everything using the theory of simplicial
categories and the notion of morphisms developped by Cordier and Porter
\cite{CordierPorter}. An alternative is to use the theory of Segal
categories,
which are a type of weak replacement for the notion of simplicial
categories. Any Segal category is equivalent to a simplicial category
\cite{DwyerKanSmith},
and one assumes that the morphisms as defined in \cite{ahcs}
are the same as the morphisms of simplicial categories defined using the
coend construction in \cite{CordierPorter} (I am not sure that anybody
has actually checked this yet). 
There is a notion of adjointness of functors between
simplicial or Segal categories. See Cordier-Porter \cite{CordierPorter}, and \cite{ahcs}.

There are actually two ways of constructing $\Ll_n$ starting with a theory of $n$-categories.
One way as mentionned above is to have a theory in which one first
constructs a $1$-category of objects denoted $nCat$, with strict
morphisms (in particular the morphisms between objects in $nCat$ are not
the real ``weak'' morphisms, which should form an $n$-category rather
than a set anyway); some of
the morphisms in $nCat$ are equivalences of $n$-categories; and define $\Ll_n$
to be the Dwyer-Kan localization of $nCat$ dividing by the weak
equivalences.
The other way is to have a theory in which there is an $n+1$-category
$nCAT$ of $n$-categories, with presumably the correct mapping
$n$-categories between objects. In this case, define $\Ll_n$ to be the
{\em $1$-groupic interior} of $nCAT$, namely the $n+1$-category with the
same objects, the same $1$-morphisms, but whose $i$-morphisms for $i\geq
2$ are only the ones which are invertible up to equivalence. Thus if $A$
and $B$ are $n$-categories,
$$
Hom_{\Ll _n}(A,B) \subset Hom_{nCAT}(A,B)
$$
is the {\em interior} or the sub-$n$-category with the same objects 
but with only invertible $i$-morphisms for $i\geq 1$.

One expects $nCat$ (or some related $1$-category) to have a closed 
model structure. If this is a simplicial closed model structure then
$\Ll _n$ will also be equivalent to the simplicial category of fibrant 
and cofibrant objects here.

For the theory of Tamsamani $n$-categories, all of the above
descriptions hold and are known to give equivalent Segal categories
which we denote by $\Ll _{n,Ta}$. 
If ${\bf T}$ is another theory of $n$-categories as in
Leinster's list \cite{Leinster} or others \cite{May} etc., then there is
always a $1$-category $nCat_{\bf T}$ and a notion of equivalence (an
equivalence is a functor which is fully faithful and essentially
surjective)
and we set $\Ll _{n,{\bf T}}$ equal to the Dwyer-Kan localization of
$nCat_{\bf T}$ divided by the equivalences.  One hopes 
in the case of other theories to have an analogue of the second
description as the interior of $nCAT_{\bf T}$ too.

The rather limited purpose of this note is to consider the categories $\Ll_n$ obtained
as above, as abstract simplicial categories, and to ask what properties
they are expected to have. Our comparison conjecture says that our list
of
properties \ref{p1}--\ref{p10}
serves to characterize $\Ll_n$ up to equivalence of
simplicial categories. We also make a conjecture about the automorphisms
of $\Ll_n$, which basically says that the equivalence in question will
be unique up to the action of $(\zz /2)^n$ obtained using the
operations of reversing the directions of the
$i$-arrows for different $i$. 

We leave as 
an exercise for the reader, to prove properties \ref{p1}--\ref{p10} in the case of
$\Ll _{n,Ta}$. For the other theories ${\bf T}$ we conjecture that 
properties \ref{p1}--\ref{p10} hold for $\Ll _{n,{\bf T}}$.

So, for the rest of the paper, we will just fix a Segal category $\Ll_n$
which is supposed to be the Segal category of weak $n$-categories
without reference to which particular theory is concerned, and
we will state the properties which it is supposed to have.

\begin{property}
\label{p1}
The category $\Ll_n$ is $n+1$-truncated, i.e. the $Hom$-spaces
between any two objects have vanishing homotopy groups in degrees $>n$.
\end{property}

\begin{property}
\label{p2}
The category $\Ll_n$ admits small homotopy limits and
colimits.
\end{property}

More generally it is an $\infty$-pretopos in the sense of \cite{Giraud}
but it isn't clear whether we need to know this.

{\em Notations:} Let $\ast$ denote the final object (which is the empty homotopy
limit)
and $\emptyset$ the initial object (which is the empty homotopy
colimit).
Denote by $k\ast$ the disjoint union of $k$ copies of $\ast$. 
If $X\in {\rm ob} \Ll _n$ then put
$$
\pi _0(X):= \pi _0 Hom (\ast , X).
$$
A morphism $f:X\rightarrow Y$ in $\Ll_n$ is {\em essentially surjective}
if $\pi _0(f)$ is surjective.

Let $k\ast / \Ll _n$ denote the Segal category of diagrams $k\ast
\rightarrow X$ in $\Ll_n$. More precisely, if $I$ denotes the Segal
category with two objects $0$ and $1$ and one arrow from $0$ to $1$,
then $k\ast / \Ll_n$ is defined to be the fiber over $k\ast$ of the map
``evaluation at $0$'',
$$
Hom (I,\Ll _n) \rightarrow Hom (\{ 0\}, \Ll _n)=\Ll_n .
$$
Let $k\ast // \Ll_n$ denote the full sub-Segal category of objects 
$k\ast \rightarrow X$ which are essentially surjective, i.e. where
$k\ast$ surjects onto $\pi _0(X)$.

\begin{property}
\label{p3}
The inclusion $k\ast // \Ll _n\rightarrow k\ast / \Ll _n$
admits a right homotopy adjoint 
$$
tr_k : k\ast / \Ll _n\rightarrow k\ast // \Ll _n .
$$
\end{property}

The case $k=2$ is the most important: $2\ast // \Ll _n$ is the Segal
category of $n$-categories with two objects; the case where the two
objects happen to be equivalent is the case of $2\ast \rightarrow X$
with $\pi _0(X) = \ast$. The notation $tr_k$ refers to the adjoint map
as ``trimming off the additional objects''.

\begin{property}
\label{p4} 
There are four nonempty saturated
\footnote{Saturated means that if
an object is in it then so is every equivalent object.}  
full sub-Segal categories 
$$
{\cal A}, \, {\cal B}, \, {\cal C}, \, {\cal D}\, \subset 2\ast // \Ll
_n
$$
and the set of
objects of $2\ast // \Ll _n$ is a disjoint union of the objects of
${\cal A}$, ${\cal B}$, ${\cal C}$, and ${\cal D}$.
This decomposition  has the following properties: 
\newline
$(dec_1)$ the only objects of ${\cal A}$ are initial objects of $2\ast // \Ll _n$;
\newline
$(dec_2)$ there are no morphisms
from objects in ${\cal C}$ or ${\cal D}$  to objects in ${\cal B}$;
\newline
$(dec_3)$ there are no morphisms
from objects in ${\cal B}$ or ${\cal D}$  to objects in ${\cal C}$;
\newline
$(dec_4)$ if $f$ is an object of ${\cal B}$ and $g$
is an object of ${\cal B}$ or ${\cal D}$ then there is a morphism from 
$f$ to $g$; and 
\newline
$(dec_5)$ if $f$ is an object of ${\cal C}$ and $g$
is an object of ${\cal C}$ or ${\cal D}$ then there is a morphism from 
$f$ to $g$.
\end{property}

Concretely for any specific theory of 
$\Ll_n$ these four subcategories are defined as follows:
\newline
${\cal A} \subset 2\ast // \Ll _n$ is the full subcategory consisting
only of the initial object;
\newline
${\cal B} \subset 2\ast // \Ll _n$ is the full subcategory consisting
of maps $\{ x,y\} \rightarrow X$ where $X$ has at least one morphism
going from $x$ to $y$ but no morphisms going from $y$ to $x$;
\newline
${\cal C} \subset 2\ast // \Ll _n$ is the full subcategory consisting
of maps $\{ x,y\} \rightarrow X$ where $X$ has at least one morphism
going from $y$ to $x$ but no morphisms going from $x$ to $y$;
\newline
${\cal D} \subset 2\ast // \Ll _n$ is the full subcategory consisting
of maps $\{ x,y\} \rightarrow X$ where $X$ has at least one morphism
going from $x$ to $y$ and at least one morphism going from $y$ to $x$.

\begin{lemma}
Given any Segal category $\Mm$, if there exists a decomposition of $\Mm$
into four saturated full subcategories ${\cal A}$, ${\cal B}$,
${\cal C}$ and ${\cal D}$ as above (such that the object set of $\Mm$ is
a disjoint union of the object sets of these categories) and if this
decomposition satisfies properties $(dec_1)\ldots (dec_5)$ then it is
the unique decomposition satisfying those properties, up to permutation
of the factors ${\cal B}$ and ${\cal C}$. 
\end{lemma}
{\em Proof:}
%%%
%%% put in here a Coq source file with the proof in a future revision
%%%
Suppose that ${\cal A}$, ${\cal B}$,
${\cal C}$ and ${\cal D}$ are one such decomposition
of $\Mm$, and suppose that ${\cal A}'$, ${\cal B}'$,
${\cal C}'$ and ${\cal D}'$ are another one. Clearly ${\cal A} = {\cal
A}'$ by the saturation condition and $(dec_1)$. Suppose that $x \in
{\cal B}$. We claim that $x\in {\cal B}' \cup {\cal C}'$. If not, then
$x\in {\cal D}'$. In that case, choosing elements 
$x' \in {\cal B}'$ and $y'\in {\cal C}'$ respectively, we have maps
$x'\rightarrow x$ and $y'\rightarrow x$ by $(dec_4)$ and $(dec_5)$.
But by $(dec_2)$ this implies that $x'\in {\cal B}$ and $y' \in {\cal
B}$.
In particular, again by $(dec_4)$ there is a map $x' \rightarrow y'$
contradicting $(dec_3)$ for $y'$. This contradiction to our assumption
that $x\in {\cal D}'$ proves that $x\in {\cal B}'\cup {\cal C}'$.
Thus ${\cal B}\subset {\cal B}'\cup {\cal C}'$. A similar argument
proves the same for ${\cal C}$ so we get
$$
{\cal B}\cup {\cal C}\subset {\cal B}'\cup {\cal C}',
$$
and running the argument in the opposite direction gives the opposite
inclusion so 
$$
{\cal B}\cup {\cal C}=  {\cal B}'\cup {\cal C}'.
$$
This implies that ${\cal D} = {\cal D}'$. 

Now fix $x_0\in {\cal B}$. By permuting the factors ${\cal B}'$ and
${\cal C}'$ if necessary (note that this preserves the hypotheses and
conclusion of the theorem) we can assume that $x_0\in {\cal B}'$. 
Now for any $y\in {\cal C}$, we claim that $y\in {\cal C}'$; for if not
then $y$ would be in ${\cal B}'$ so there would be a map $y\rightarrow
x_0$ by $(dec_4)$, contradicting $(dec_2)$; this proves the claim. 
Thus ${\cal C}\subset {\cal C}'$. The same argument in the other
direction shows the reverse inclusion so ${\cal C}={\cal C}'$ which in
turn implies that ${\cal B} = {\cal B}'$. 
\eop

In view of this lemma, the existence
of the decomposition is a property, and up to permutation of the factors
the full subcategory ${\cal B}$ is uniquely determined (to be precise,
there are two full subcategories which can play this role). 

\begin{property}
\label{p5}
In the case of the decomposition of $2\ast // \Ll _n$ there  is
an involution of $\Ll _n$ taking a category to its opposite, which
permutes the two factors in the above decomposition. 
\end{property} 

Thus, without loss
of generality we may isolate the factor ${\cal B}$ in the decomposition
as being well defined.

Let $r_0$ and $r_1$ denote the two maps from 
$2\ast / \Ll _n$ to $\ast / \Ll _n$ obtained by retaining only the first
or second point. These restrict to maps denoted the same way on ${\cal
B}$  and we can compose with the ``trimming'' map $tr_1$ to obtain
$$
tr_1 r_i : {\cal B} \rightarrow \ast // \Ll _n .
$$
Let ${\cal B}^e$ denote the homotopy fiber of 
$$
(tr_1r_0, tr_1r_1):{\cal B} \rightarrow 
\ast // \Ll _n \times \ast // \Ll _n 
$$
over the object $(\ast , \ast )$. Note that this object (plus those
equivalent to it) constitutes a full sub-Segal category of 
$\ast // \Ll _n \times \ast // \Ll _n $ which is trivial i.e. all
mapping spaces are reduced to $\ast$, so the homotopy fiber ${\cal B}^e$
is also a full sub-Segal category of ${\cal B}$. 

Concretely in the examples, ${\cal B}^e$ is the
category of $n$-categories with two objects $0,1$ and some morphisms between
them (going only in one direction, say from $0$ to $1$), 
but where the trimmings to either endpoint are trivial. Heuristically,
this last condition means that the mapping $n-1$-categories of
endomorphisms of either $0$ or $1$ are trivial, so such an object is
determined entirely by the $n-1$-category of morphisms from $0$ to $1$.
The following property makes this precise.

Let ${\cal A} + {\cal B}^e$ denote the saturated full subcategory of
$2\ast // \Ll _n$ whose object set is the union of the object sets of
${\cal A}$ and ${\cal B}^e$ (i.e. we add to ${\cal B}^e$ the initial
object
which corresponds in the examples
to cases where the morphism $n-1$-category from $0$ to $1$ is empty).

\begin{property}
\label{p6} 
The Segal category ${\cal A}+{\cal B}^e$ is equivalent to $\Ll _{n-1}$.
\end{property} 

Choose an equivalence between ${\cal A}+{\cal B}^e$ and $\Ll _{n-1}$. 
Denote by 
$$
\Upsilon : \Ll _{n-1}\rightarrow 2\ast // \Ll _n \rightarrow
2\ast / \Ll _n
$$ 
the composed functor. 
\footnote{In the case $\Ll _n = \Ll _{n,Ta}$ it is possible to make the
choice of equivalence in Property \ref{p6}
so that the functor $\Upsilon$ is just the localization of the
functor of the same name used for example in \cite{limits}.}

\begin{conjecture} 
\label{auts}
The Segal groupoid of auto-equivalences $Aut(\Ll _n)$
is equal to $(\zz /2)^n$ with the first factor being the involution 
referred to in Property \ref{p5} above and with the remaining factors restricting
via $\Upsilon$ to
the automorphisms of $\Ll _{n-1}$.
\end{conjecture}

This is stated as a conjecture because even for $\Ll _{n,Ta}$ it is a
more advanced exercise. Nevertheless, it explains why the problem of
choice of an equivalence in Property \ref{p6} is  not too serious.

\begin{property}
\label{p7}
The functor $\Upsilon$ has a right adjoint
$$
H: 2\ast / \Ll _n \rightarrow \Ll _{n-1}.
$$
\end{property}

This is to be viewed as the functor which to a diagram $\{ x,y\}
\rightarrow X$ associates the $n-1$-category of morphisms from $x$ to
$y$ in $X$. In cases where it makes sense to talk about this
$n-1$-category of morphisms, it is a property that $H$  has this
interpretation. \footnote{That is to say, in $\Ll _{n,Ta}$ with the
choice referred to in the previous footnote, this interpretation is
true; and for any theory ${\bf T}$ in which we have an $n-1$-category of
morphisms, it is a conjecture that for an appropriate choice of
equivalence in Property \ref{p6} above, this interpretation is true.}

\begin{property}
\label{p8}
The adjunction morphism for $E \in \Ll _{n-1}$ is an
equivalence
$$
E \stackrel{\cong}{\rightarrow} H (\Upsilon (E)).
$$
\end{property}

Say that a morphism $f:X\rightarrow Y$ in $\Ll _n$ is {\em fully faithful}
if for every $u:\{ x,y\} \rightarrow X$ the resulting map
$$
H(X,u)\rightarrow H(Y,f\circ u)
$$
is an equivalence in $\Ll _{n-1}$.

This notion coincides with the previously existing notions of fully
faithful map in $\Ll_n$. 

\begin{property}
\label{p9}
Equivalence in $\Ll _n$ is characterized using this notion: a morphism $f$ is an
equivalence if and only if it is essentially surjective and fully
faithful in the above senses.
\end{property}

In particular, this means that an object $X$ with $\pi _0(X)=\emptyset$
is equivalent to the initial object $\emptyset$. It is interesting to
note that the Segal category $\Mm$ of locally constant $n$-stacks on the
$n+1$-sphere satisfies most of the other properties we write down here,
but there are $n$-stacks classified by nonzero cohomology classes in
degree $n+1$ which have no global sections, thus have $\pi _0=\emptyset$
but are not equivalent to $\emptyset$. Property \ref{p9} is an important
way of ruling out that type of example.

Given $E_1,\ldots , E_k\in \Ll _{n-1}$, define the homotopy coproduct
$$
\Upsilon ^k(E_1,\ldots , E_k) := 
\Upsilon (E_1) \cup ^{\{ 1\}} \ldots \cup ^{\{ k-1\}} \Upsilon (E_k)
$$
where we have denoted the $k+1$ objects $\ast$ by $0,\ldots , k$.
These are organized as follows: the first and last objects of the term
$\Upsilon (E_i)$ are $i-1$ and $i$. Thus this coproduct represents 
``composing arrows head to tail''.
The main case is 
$$
\Upsilon ^2(E,F) = \Upsilon (E) \cup ^{\{ 1\}} \Upsilon (F).
$$

The morphism $\Upsilon (E_i)\rightarrow \ast$ induces a morphism of the
coproduct
$$
\Upsilon ^k(E_1,\ldots , E_k) \rightarrow \Upsilon ^{k-1}(E_1,\ldots ,
\widehat{E_i},
\ldots , E_k).
$$
In particular, composing many of these we get down to morphisms
$$
\Upsilon ^k(E_1,\ldots , E_k) \rightarrow \Upsilon (E_j)
$$
for any $1\leq j \leq k$.
\footnote{In the case of $\Ll _{n,Ta}$ with the good choice made in 
Property \ref{p6},
the construction $\Upsilon ^k$
and the various maps are induced by  those considered in \cite{limits}.}
Given that $H(\Upsilon (E))\cong E$ by Property \ref{p8}, we obtain a morphism
$$
H(\Upsilon ^k(E_1,\ldots , E_k), \{ 0,k\} ) \rightarrow
E_1 \times \ldots \times E_k
$$
where $\times$ denotes the direct product in $\Ll _{n-1}$.

\begin{property}
\label{p10}
{\sc (main property)} The above morphism is an equivalence: 
$$
H(\Upsilon ^k(E_1,\ldots , E_k), \{ 0,k\} ) \stackrel{\cong}{\rightarrow}
E_1 \times \ldots \times E_k .
$$
\end{property}

This finishes our list of properties.

\begin{conjecture}
\label{char}
The above properties \ref{p1}--\ref{p10} characterize the Segal
category $\Ll _n$ up to equivalence (the choice of equivalence is
governed by Conjecture \ref{auts}). 
\end{conjecture}

This conjecture would give a way of approaching the comparison problem:
to compare two definitions of $n$-category, if we could establish
properties \ref{p1}--\ref{p10} for both of them, and if we knew the
conjecture, then we would get equivalence of the definitions at least in
the sense of equivalence between the Dwyer-Kan localized homotopy categories.

To finish this note, we just show how the above properties lead to the
composition map 
$$
H(X,\{ x , y\} ) \times H (X,\{ y, z\}) \rightarrow H(X,\{ x, z\}).
$$
Suppose $X$ is an object of $\Ll _n$ and suppose 
$x,y,z: \ast \rightarrow X$ are three maps. Let 
$(X,\{ x,y\})$ denote the object of $2\ast / \Ll _n$ corresponding to
$(x,y): 2\ast \rightarrow X$ and so forth. Put
$$
E:= H(X,\{ x, y\} ), \;\;\; F := H(X,\{ y,z\} ).
$$
By adjunction we obtain maps $\Upsilon (E)\rightarrow X$ and 
$\Upsilon (F)\rightarrow X$ sending $0,1,2$ to $x,y,z$ respectively
(here the two objects of $\Upsilon (F)$ are denoted $1,2$ as per the
convention above). These induce a map on the coproduct 
$$
\Upsilon ^2(E,F) = \Upsilon (E) \cup ^{\{ 1\} }\Upsilon (F) 
\stackrel{\alpha}{\rightarrow}
X.
$$
On the other hand, Property \ref{p10} says that the map 
$$
H(\Upsilon ^2(E,F), \{ 0,2\} )\stackrel{\beta}{\rightarrow} E \times F
$$
induced by the maps $\Upsilon ^2 (E,F)\rightarrow \Upsilon (E)$ and
$\Upsilon ^2(E,F)\rightarrow \Upsilon (F)$, is an equivalence. 
The map $\alpha$ induces a map
$$
H(\Upsilon ^2(E,F), \{ 0,2\} )\rightarrow 
H(X,\{ x, z\} ),
$$
and composing with the inverse of the equivalence $\beta$ we obtain a
map
$$
E\times F \rightarrow H(X,\{ x,y\} ).
$$
Recalling the definitions of $E$ and $F$ this becomes a map
$$
\gamma (x,y,z):H(X,\{ x , y\} ) \times H (X,\{ y, z\}) \rightarrow H(X,\{ x,
z\} ).
$$
This map is the ``composition'' map between the mapping objects. The
definition of the mapping object functor $H$ and this composition have
been constructed using only the structure of the simplicial category
$\Ll _n$. 

We expect that, after solving the homotopy-coherence problems involved,
one would obtain more generally a functor from objects of $\Ll _n$ to
simplicial objects in $\Ll _{n-1}$ having discrete $0$-th stage and
satisfying the Segal properties; this would allow one to obtain
(by induction) a
functor 
$$
\Ll _n \rightarrow \Ll _{n,Ta}
$$
and in order to prove Conjecture \ref{char} it would suffice
to prove that this map is an equivalence.

\end{document}